\documentclass[a4paper,reqno]{article}
\usepackage{amsmath}
\usepackage{amssymb}
\usepackage{amsthm}
\newcommand{\g}{\mathfrak{g}}

\newcommand{\tb}{\mathbf{t}}
\newcommand{\Uh}{{U}_h (\g)}

\newcommand{\A}{\mathcal{A}}

\newcommand{\C}{\mathbb{C}}

\newcommand{\Z}{\mathcal{Z}}
\newcommand{\X}{\mathcal{X}}
\newcommand{\I}{\mathcal{I}}
\newcommand{\ot}{\otimes}

\newcommand{\ad}{\operatorname{ad}}
\newcommand{\id}{\operatorname{id}}
\newcommand{\Tr}{\operatorname{Tr}}
\newcommand{\End}{\operatorname{End}}
\newcommand{\bL}{\overline{L}}
\newcommand{\vep}{\varepsilon}
\renewcommand{\Im}{\operatorname{Im}}
\newcommand{\Ker}{\operatorname{Ker}}
\newcommand{\bsig}{\bar{\sigma}_h}
\theoremstyle{plain} 

\theoremstyle{definition}

\theoremstyle{remark}
\newtheorem{rem}{Remark}[section]

\numberwithin{equation}{section}

\newcommand{\be}[1]{\begin{eqnarray#1}}
\newcommand{\ee}[1]{\end{eqnarray#1}}

\begin{document}
\begin{center}
{\large\bf DOUBLE QUANTIZATION ON THE COADJOINT\\ 
REPRESENTATION OF $\mathbf{ SL(n)}$}
\footnote[1]{Presented at the 6th International Colloquium on Quantum Groups:
``Quantum Groups and Integrable Systems'', Prague, 19-21 June 1997.}\\
\bigskip
 J. DONIN\\
\medskip
{\it Department of Mathematics, Bar-Ilan University,\\
52900 Ramat-Gan, Israel\\
\medskip
E-mail: donin@macs.biu.ac.il}
\end{center}
\begin{abstract}
For $\g=sl(n)$ we construct a two parametric
$U_h(\g)$-invariant family of algebras, $(S\g)_{t,h}$, which defines
a quantization of the function algebra $S\g$ 
on the coadjoint representation and in the parameter $t$ 
gives a quantization of the Lie bracket. The family induces a 
two parametric deformation of the function algebra 
of any maximal orbit
which is a quantization of the Kirillov-Kostant-Souriau bracket
in the parameter $t$.
In addition we construct a quantum de Rham complex on~$\g^*$.
\end{abstract}
\section{Introduction}
Let $G$ be a simple Lie group over the field of complex numbers $\C$
with the Lie algebra $\g$.
Let $M$ be a homogeneous space and $\A$ the algebra of algebraic
functions on $M$. The algebra $\A$ is commutative and has a $G$-invariant
multiplication. 

The first problem we consider in this note is to construct a quantized 
algebra $\A_h$
in which the deformed multiplication $m_h$ is invariant under the 
action of the Drinfeld-Jimbo
quantum group $U_h(\g)$ defined over $\C[[h]]$, i.e.
\be{*} \label{inv2}
xm_h(a\ot b)=m_h\Delta_h(x)(a\ot b),
\ee{*} 
where $a,b\in\A_h$, $x\in U_h\g$, and
$\Delta_h$ denote the comultiplication in $U_h(\g)$.

Let $\A_t$ be a $G$-invariant quantization of $\A$ by an invariant
Poisson bracket on $M$. 
The second problem is to construct a two parametrized
$\Uh$-invariant family of algebras, $\A_{t,h}$, such that $\A_t=\A_{t,0}$.

More explicitly,
by a quantization (or deformation) of $\A$ we mean an algebra $\A_h$ over 
the algebra
of formal power series $\C[[h]]$ in a variable $h$
that is isomorphic to $\A[[h]]=\A\ot\C[[h]]$ as a 
$\C[[h]]$-module and
$\A_0=\A_h/h\A_h=\A$ as an algebra
(note that when we consider modules over  $\C[[h]]$
the symbol $\ot$ denotes the tensor product 
completed in the $h$-adic topology).  We say in this case
that $\A_h$ is a flat deformation of $\A$.

Similarly, the flatness of two parametric family, $\A_{t,h}$,
is defined. 

In \cite{DGM}, \cite{DG1}, \cite{DG3} the first problem is solved for  
flag varieties,
for orbits of highest weight vector in irreducible representations
of $G$, and for semisimple orbits in $\g^*$. 
The second problem is solved in \cite{DS} for hermitian 
symmetric spaces. 

In this note we consider the coadjoint representation $\g^*$
of $U(\g)$ for $\g=sl(n)$ and show that the second problem can be solved when 
$M$ is a maximal orbit in $\g^*$ and $\A_t$ is a quantization of
the Kirillov-Kostant-Souriau bracket.

Moreover, we define a $\Uh$-invariant algebra $(S\g)_{t,h}$, quantization
of the function algebra on $\g^*$, and show that the quantized algebra
$\A_{t,h}$ for a maximal orbit can be presented as a quotient algebra
of $(S\g)_{t,h}$.

To obtain $(S\g)_{t,h}$ we use an idea of the paper of Lyubashenko
and Sudbery \cite{LS} concerning
the construction of
a quantum analog of Lie algebra for $sl(n)$,
but we deal with the quantum group $\Uh$ over $\C[[h]]$ instead of
$U_q(\g)$ and use for our construction the R-matrix, that
allows us to simplify the proofs. Using the fact that
the quantum Casimir $C_V$ is invertible in $\Uh$ (see Section 2),
we define a quantum Lie algebra as an embedding of the deformed adjoint 
representation, $\g_h$, in $\Uh$ such that the kernel of the extension
to an algebra homomorphism $T(\g_h)\to\Uh$ is defined 
by quadratic-linear relations,
as in the classical case. This definition differs from the definition given in
\cite{LS}. For $U_q(\g)$ our results are valid as well (see Remark 3.2).

I am grateful to D.Gurevich and S.Shnider for helpful discussions.

\section{Quantum Lie algebra for $U_h(sl(n))$}
Let $R=R^\prime_i\ot R^{\prime\prime}_i\in \Uh\ot \Uh$ be the R-matrix 
(summation by $i$ is assumed).
It satisfies the properties \cite{Dr}
\be{} \label{Delo}
\Delta^\prime(x)=R\Delta(x)R^{-1}, \quad x\in \Uh,
\ee{}
where $\Delta$ is the comultiplication in $\Uh$ and $\Delta^\prime$  
the opposite one,
\be{} \label{six}
(\Delta\ot 1)R&=&R^{13}R^{23}=R^\prime_i\ot R^\prime_j\ot R^{\prime\prime}_i
R^{\prime\prime}_j \notag\\
(1\ot\Delta)R&=&R^{13}R^{12}=R^{\prime}_iR^{\prime}_j\ot R^{\prime\prime}_j
\ot R^{\prime\prime}_i, 
\ee{}
and
\be{} \label{coun}
(1\ot \varepsilon)R=(\varepsilon\ot 1)R=1\ot 1, 
\ee{}
where $\varepsilon$ is the counit in $\Uh$.

Consider the element $Q=Q_i^\prime\ot Q_i^{\prime\prime}=R^{21}R$. 
It follows from (\ref{Delo}) that $Q$
commutes with elements from $\Uh\ot \Uh$ of the form $\Delta(x)$. This is
equivalent to $Q$ being invariant under the adjoint action of $\Uh$
on  $\Uh\ot \Uh$.

Let $V$ be an irreducible finite dimensional representation of $\Uh$ and 
$\rho:\Uh\to\End(V)$ the corresponding map of algebras.
Consider the dual space $\End(V)^*$ as a left $\Uh$-module setting
$$(x\varphi)(a)=\varphi(\gamma(x_{(1)})ax_{(2)}), $$
where $\varphi\in\End(V)^*$, $a\in\End(V)$, $\Delta_h(x)=x_{(1)}\ot x_{(2)}$
in Sweedler notions, and $\gamma$ denotes the antipode in $\Uh$.

Consider the map $f:\End(V)^*\to\Uh$ defined as
$\varphi\mapsto \varphi(\rho(Q^{\prime}_i)Q^{\prime\prime}_i$. 
From the invariance of
$Q$ it follows that $f$ is a $\Uh$-equivariant map, so $\bL=\Im(f)$ is a
$\Uh$-submodule. 

It follows from (\ref{six}) that $\bL$ is a left coideal in $\Uh$, i.e.
$\Delta(x)\in \Uh\ot \bL$ for any $x\in \bL$. 
Indeed, $Q=R^{\prime\prime}_iR^{\prime}_j\ot R^{\prime}_iR^{\prime\prime}_j$.
Applying  
(\ref{six}) we obtain
\be{*}
(1\ot\Delta_h)R^{21}R=R^{\prime\prime}_i
R^{\prime\prime}_jR^{\prime}_kR^{\prime}_l\ot R^{\prime}_iR^{\prime\prime}_l\ot
R^{\prime}_jR^{\prime\prime}_k
\ee{*}
Let $\varphi\in \End(V)^*$. Define $\psi_{il}\in\End(V)^*$ setting
$\psi_{il}(a)=\varphi(R^{\prime\prime}_iaR^{\prime}_l)$ for $a\in\End(V)$.
Then
$\Delta\varphi(R^{\prime\prime}_iR^{\prime}_j)R^{\prime}_iR^{\prime\prime}_j=
R^{\prime}_iR^{\prime\prime}_l\ot \psi_{il}(R^{\prime\prime}_j
R^{\prime}_k)R^{\prime}_jR^{\prime\prime}_k,$
which obviously belongs to $\Uh\ot\bL$. 
\smallskip

Recall \cite{Dr}, that $R=F^{21}_he^{\frac{h}{2}\tb}F^{-1}_h$. 
Here $\tb=\sum_i t_i\ot t_i$ is the split Casimir, where 
$t_i$ form an orthonormal
basis in $\g$ with respect to the Killing form, 
$F=1\ot 1+\frac{h}{2}r+o(h)$,
and $r$ is the classical Drinfeld-Jimbo R-matrix. 
Therefore, 
\be{} \label{rmt}
Q=R^{21}R=Fe^{h\tb}F^{-1}=1\ot 1+h\tb+\frac{h^2}{2}(\tb^2+[r,\tb])+o(h^2).
\ee{}

Denote by $\Tr$ the unique (up to a factor) invariant element in $\End(V)^*$.
Let $Z_0=\rho_0(\g)$ and denote by $Z_h$ some $\Uh$-invariant deformation of
$Z_0$ in $\End(V)$.  
Then we have a decomposition $\End(V)=I\oplus Z_h\oplus W$, where $I$ is 
one dimensional subspace of invariant elements, $W$ is a subspace complement 
to $I\oplus Z_h$.
This gives a decomposition $\End(V)^*=I^*\oplus Z_h^*\oplus W^*$ where
$W^*$ consists of all the elements which are equal to zero on $I\oplus Z_h$.  
The space $I^*$ is generated by $\Tr$, and after normalizing in such a way
that $\Tr(\id)=1$, we obtain that $C_V=f(\Tr)$ is of the form
\be{} \label{rcv}
C_V=\Tr(\rho(Q_1))Q_2=1+h^2c+o(h^2),
\ee{}
and $c$ has to be an invariant element in $U(\g)$. It follows from (\ref{coun})
that $\varepsilon(C)=1$.

From (\ref{rmt}) follows that the elements of $f(Z_h^*)$ have the form
\be{}
z=hx+o(h), \quad x\in \g,
\ee{}
hence the subspace $L_1=h^{-1}f(Z_h^*)$ forms a subrepresentation in $\Uh$
under the left adjoint action of $\Uh$ on itself, which is a deformation
of the standard embedding of $\g$ into $U(\g)$. It follows from (\ref{coun})
that $\varepsilon(L_1)=0$.

The elements from $f(W^*)$ have the form $w=h^2b+o(h^2)$
and $\varepsilon(W^*)=0$. 

Denote $L=h^{-1}f(Z_h^*+W^*)$, so 
$\bL=\C C_V\oplus hL$.
Since $\bL$ is a left coideal in $\Uh$, for any $x\in L$ we have
\be{*} 
\Delta(x)=x_{(1)}\ot x_{(2)}=z\ot C_V+v\ot x^{\prime},
\ee{*}
where $z,v\in \Uh$, $x^{\prime}\in L$.
Applying to the both sides $(1\ot\varepsilon)$ and multiplying
we obtain $x=x_{(1)}\vep(x_{(2)})=z\vep(C_V)+v\vep(x^{\prime})=z$. So, $z$ have
to be equal to $x$. and we obtain
\be{} \label{qlr1}
\Delta(x)=x_{(1)}\ot x_{(2)}=x\ot C_V+v\ot x^{\prime}, 
\quad x,x^{\prime}\in L.
\ee{}

From (\ref{qlr1}) we have for any $y\in L$ 
\be{} \label{qlr2}
xy=x_{(1)}y\gamma(x_{(2)})x_{(3)}=x_{(1)}y\gamma(x_{(2)})C_V+v_{(1)}y
\gamma(v_{(2)})x^{\prime}.
\ee{}

Introduce the following maps:
\be{} \label{qlr'}
\sigma_h^{\prime}:&L\ot L\to L\ot L,&\quad x\ot y\mapsto 
v_{(1)}y\gamma(v_{(2)})\ot x^{\prime},  \notag \\  
\ [\cdot,\cdot]_h^{\prime}:&L\ot L\to L, &\quad x\ot y\mapsto x_{(1)}y
\gamma(x_{(2)}).
\ee{}
We may rewrite (\ref{qlr2}) in the form
\be{} \label{qlr3}
m(x\ot y-\sigma_h^{\prime}(x\ot y))-[x,y]_h^{\prime}C_V=0.
\ee{}

Observe now that it follows from (\ref{rcv}) that 
$C_V$ is an invertible element in $\Uh$. Put $P=C^{-1}_V$.
Transfer the maps (\ref{qlr'}) to the space $P\cdot L$, i.e. define
\be{*}
\sigma_h(Px,Py)&=&(P\ot P)\sigma^{\prime}_h(x,y), \\
{[}Px,Py]_{h} &=&P[x,y]^{\prime}_h.
\ee{*}
From (\ref{qlr1}) we obtain
\be{}
P_{(1)}x_{(1)}\ot P_{(2)}x_{(2)})=P_{(1)}x\ot P_{(2)}C_V+P_{(1)}v\ot 
P_{(2)}x^{\prime},
\ee{}
Using this relation and taking into account that $P$ commutes with all 
elements from $\Uh$, we obtain as in (\ref{qlr2}) 
\be{*}
PxPy&=&P_{(1)}x_{(1)}Py\gamma(x_{(2)})\gamma(P_{(2)})P_{(3)}x_{(3)}= \\
& &P_{(1)}x_{(1)}Py\gamma(x_{(2)})\gamma(P_{(2)})P_{(3)}C_V+P_{(1)}v_{(1)}Py 
\gamma(v_{(2)})\gamma(P_{(2)})P_{(3)}x^{\prime}= \\
& & P[x,y]^{\prime}_h+P^2m\sigma_h^{\prime}(x\ot y)= 
 [Px,Py]_h+m\sigma_h(Px\ot Py)). 
\ee{*}
This equality may be written as
\be{} \label{qlr4}
m(x\ot y-\sigma_h(x\ot y))-[x,y]_h=0, \quad x,y\in C^{-1}_VL. 
\ee{}

Define $L_V=C^{-1}_VL$.
Let $T(L_V)=\oplus_{k=0}^\infty L_V^{\ot k}$ be the tensor algebra over $L_V$.
Notice, that $T(L_V)$ is not supposed to be completed in $h$-adic topology.
Let $J$ be the ideal in $T(L_V)$ generated by the relations
\be{} \label{rel0}
(x\ot y-\sigma_h(x\ot y))-[x,y]_h, \quad x,y\in L_V.
\ee{}
Due to (\ref{qlr4}) we have a homomorphism of algebras,
$T(L_V)/J\to \Uh$, extending the natural embedding $\imath:L_V\to \Uh$.
Introduce a new variable $t$ and consider a homomorphism of algebras,
$T(L_V)[t]\to \Uh[t]$, which extends the embedding $t\imath:L_V[t]\to \Uh[t]$.
From (\ref{qlr4}) follows that it factors through the homomorphism
of algebras
\be{} \label{phit}
\phi_{t,h}:T(L_V)[t]/J_t\to \Uh[t],
\ee{} 
where $J_t$ is the ideal generated by the relations
\be{} \label{relt}
(x\ot y-\sigma_h(x\ot y))-t[x,y]_h, \quad x,y\in L_V.
\ee{}

Let now $\g=sl(n)$, and $V$ be a basic (defining) representation of $\g$.
In this case $L_V$ is isomorphic to a deformed standard embedding of
$\g$ into $U(\g)$. In the next section we shall see that the quantum group
$\Uh$ as an algebra is defined via this embedding by some  
$\Uh$-invariant quadratic-linear relations in
$T(L_V)$,
as in the classical case. For this reason we call $L_V$ a quantum Lie
(sub)algebra of $U_h(sl(n))$.

\section{Double quantization on $sl(n)^*$ and quantum de Rham complex} 
In this section $\g=sl(n)$.

Let us apply the construction of the previous section to $V=\C^n[[h]]$,
the deformed basic representation of $\g$. In this case $\End(V)=I\oplus Z_h$, 
where $Z_h$ is a deformed adjoint representation.
So, $\g_h=L_V=h^{-1}C_V^{-1}f(Z^*_h)$ is a deformation of the standard 
embedding of $\g$ in $U(\g)$. It is easy to see that in this case
$\sigma_h$ is a deformation of the usual permutation: 
$\sigma_0(x\ot y)=y\ot x$, and $[\cdot,\cdot]_h$ is a deformation of the
Lie bracket on $\g$: $[x,y]_0=[x,y]$, $x,y\in \g\subset U(\g)$.

It follows that the homomorphism of algebras
$T(\g_h)/J\to \Uh$ is a monomorphism, because it is an isomorphism for $h=0$.
Recall, that the ideal $J$ is defined by relations (\ref{rel0}).
Moreover, according to the PBW theorem the algebra $\Im(\phi_{t,h})$ 
at the point $h=0$ is a free $\C[t]$-module and is equal to
\be{*} \label{qst}
(S\g)_t=T(\g)/\{x\ot y-y\ot x-t[x,y]\}. 
\ee{*}
For $t=0$ this algebra is the
symmetric algebra $S\g$, the algebra of algebraic functions on $\g^*$. 
For $t\neq 0$ this algebra 
is isomorphic to $U(\g)$. Moreover, $\Uh$ is a flat $\C[[h]]$-module.
It follows 
that $\phi_{t,h}$ in (\ref{phit}) is a monomorphism of algebras over 
$\C[[h]][t]$
and $\Im(\phi_{t,h})$ is a free $\C[[h]][t]$-module isomorphic to
\be{} \label{qsth}
(S\g)_{t,h}=T(\g_h)[t]/\{x\ot y-\sigma_h(x\ot y)-t[x,y]_h\}. 
\ee{}

Call the algebra 
\be{} \label{qsh}
(S\g)_h=(S\g)_{0,h}=T(\g_h)/\{x\ot y-\sigma_h(x\ot y)\}
\ee{}
a quantum symmetric algebra. It is a 
free $\C[[h]]$-module and a quadratic algebra equal to $S\g$ at $h=0$.
 
\begin{rem} 
The completion of $(S\g)_{1,h}$ in $h$-adic topology is isomorphic to
the quantum group.
Hence, the quantum group $\Uh$ for $\g=sl(n)$  may by
considered in some sense as a quadratic-linear algebra.

Two compatible Poisson brackets correspond to the deformation $(S\g)_{t,h}$.
One of them is a quadratic one defined by the operator
$\sigma_h$. 
Another one is the usual Lie bracket.
\end{rem}

Now define a quantum exterior algebra, $(\Lambda\g)_h$.

First, modify the operator $\sigma_h$. 
Since the representation $\g^*_h$ is 
isomorphic to $\g_h$, there exists a $\Uh$-invariant bilinear form on $\g_h$,
deformed Killing form. 
This form can be extended to all tensor degrees $\g^{\ot k}$.
Let $W_h^2$ be the $\C[[h]]$-submodule in $\g_h\ot \g_h$
orthogonal to $V_h^2=\Im(\id\ot\id-\sigma_h)$. 
Define an operator $\bsig$ on $\g_h\ot\g_h$ in such a way that it has
the eigenvalues $-1$ on $V_h^2$ and
$1$ on $W_h^2$. It is clear, that $V_h^2$ and $W_h^2$ are 
deformed skew symmetric
and symmetric subspaces of $\g\ot \g$. 

Now observe, that the third graded component in the quadratic algebra
$(S\g)_h$ is the quotient of $\g_h^{\ot 3}$ by the submodule 
$V_h^2\ot\g_h+\g_h\ot V_h^2$, hence this submodule and, 
therefore, the submodule
$V_h^2\ot\g_h\cap \g_h\ot V_h^2$ are direct summands in $\g_h^{\ot 3}$, i.e. 
they have
complement submodules. As the complement submodules one can
choose the submodules $W_h^2\ot\g_h\cap \g_h\ot W_h^2$ and 
$W_h^2\ot\g_h+\g_h\ot W_h^2$, respectively, since they are complement 
at the point $h=0$
and $W_h^2$ is dual to $V_h^2$ with respect to the Killing form extended
to $\g_h\ot\g_h$. 
Hence, $W_h^2\ot\g_h+\g_h\ot W_h^2$ is a direct submodule. 
Moreover, the symmetric algebra $S\g$ is Koszul. From
a result of Drinfeld \cite{Dr2} it follows that the quadratic algebra 
$(\Lambda\g)_h=T(\g_h)/\{W_h^2\}$ is a free $\C[[h]]$-module, i.e. is
a flat $\Uh$-invariant deformation of the exterior algebra $\Lambda\g$.

Call $(\Lambda\g)_h$ a quantum exterior algebra over $\g$.

Define a quantum algebra of differential forms over $\g^*$ as the tensor
product $(\Omega\g)_h=(S\g)_h\ot(\Lambda\g)_h$ in the tensor category of
representations of the quantum group $\Uh$. The multiplication of two
elements $a\ot\alpha$ and $b\ot\beta$ looks like $ab_1\ot\alpha_1\beta$,
where $b_1\ot\alpha_1=S(\alpha\ot b)$ and $S=\sigma R$ is the permutation
in that category.

As in the classical case, the algebras $(S\g)_h$ and $(\Lambda\g)_h$
can be embedded in $T(\g_h)$ as a graded submodules
in the following way. Call the submodule 
$W^k_h=(W_h^2\ot\g_h\ot\cdots\ot\g_h)\cap
(\g_h\ot W_h^2\ot\g_h\ot\cdots\ot\g_h)\cap\cdots\cap
(\g_h\ot\g_h\ot\cdots\ot W_h^2)$
of $T^k(\g_h)$ a $k$-th symmetric part of $T(\g_h)$.
It is clear that the natural map $\pi_W:T(\g_h)\to (S\g)_h$ restricted
to $W^k_h$ is a bijection onto the $k$-degree component $(S^k\g)_h$ 
of $(S\g)_h$. Denote by $\pi_W^{\prime}:(S^k\g)_h\to W^k_h$
the inverse bijection. Similarly we define $V^k_h$, the $k$-th 
skew symmetric part of $T(\g_h)$, and the bijection
$\pi_V^{\prime}:(\Lambda^k\g)_h\to V^k_h$.

Now, define a differential $d_h$ in $(\Omega\g)_h$ as a homogeneous
operator of degree $(-1,1)$. It acts on the element $a\ot\omega$
of degree $(k,m)$ in the following way. Let $a\ot\omega=
(a_1\ot\cdots\ot a_k)\ot(\omega_1\ot\cdots\ot\omega_m)$ be
its realization as an element from $W^k_h\ot V^m_h$.
Then the formula 
\be{}
d_h(a\ot\omega)=(a_1\ot\cdots\ot a_{k-1}\ot\pi_V^\prime\pi_V
(a_k\ot\omega_1\ot\cdots\ot\omega_m)
\ee{}
presents the element $d_h(a\ot\omega)$ through its realization in
$W^{k-1}_h\ot V^{m+1}_h$. One can prove that $d_h^2=0$.

Call the algebra $(\Omega\g)_h$ with operator $d_h$ a 
quantum de Rham complex. 
It is easy to see that at the point $h=0$ this complex becomes 
the usual de Rham complex. 
The quantum de Rham complex is exact, since $d^2_h=0$ and
it is exact at $h=0$.

\begin{rem} Up to now all our constructions were considered for
the quantum group in sense of Drinfeld, $\Uh$, defined over $\C[[h]]$.
But one can deduce all the constructions above for
the quantum group in sense of Lusztig, $U_q(\g)$, defined over
the algebra $\C[q,q^{-1}]$.
We show, for example, how to obtain the quantum
symmetric algebra over $\g$.
Let $E$ be a Grassmannian consisting of subspaces  
in $\g\ot\g$ of dimension equal to $\dim(\Lambda^2\g)$,
and $\Z$ the closed algebraic subset of $E$ consisting of 
subspaces $J$ such that
$\dim(E\ot J\cap J\ot E)\geq\dim(\Lambda^3\g)$.
Let $\X$ be 
the algebraic subset in $\Z\times(\C\setminus 0)$   
consisting of points $(J,q)$ such that $J$
is invariant under the action of $U_q(\g)$. The projection
$\X\to(\C\setminus 0)$ is a proper map.
One can check that over the point $q=1$ there lies only one point of $\X$. 

As follows from the existence of $(S\g)_h$ (completed situation at $q=1$),
the dimension of $\X$ is equal to $1$. Hence, the projection 
$p:\X \to\C\setminus 0$ is a covering. For $x\in\X$ let $J_x$ be the
corresponding subspace in $\g\ot\g$ and $(S\g)_x=T(\g)/\{J_x\}$
the corresponding quadratic algebra. Due to the projection $p$ the
family $(S\g)_x$, $x\in\X$, is a module over $\C[q,q^{-1}]$. Since
$J_x$ is $U_{p(x)}(\g)$-invariant, $(S\g)_x$ is a $U_{p(x)}(\g)$-algebra.
At the ``classical'' point $x_0$, $p(x_0)=1$, this module is flat. 
Hence, after possibly deleting from
$\X$ some countable
set of points, we obtain a flat family of algebras with the same
Poincar\'{e} series as $S\g$. So, $(S\g)_x$ is the quantum symmmetric
algebra over $U_q(\g)$.
\end{rem}
\section{Double quantization on semisimple orbits in $sl(n)^*$}
In this section $G=SL(n)$, $\g=sl(n)$.

Let $M$ be a semisimple orbit of $G$ in $\g^*$ and $A$ the algebra of algebraic
functions on $M$. It is known that $M$ is a closed algebraic submanifold
in $\g^*$ \cite{Dix}, so $A$ can be presented as a quotient of $S\g$ by some 
ideal,
$S\g\to A\to 0$. The Lie bracket on $\g^*$ induces the Kirillov-Kostant-Souriau
(KKS) bracket on $M$. 

The problem is to construct a two parametric family of algebras, $A_{t,h}$,
such that
$A_{t,h}$ is a free $\C[[h]][t]$-module, $A_{0,0}=A$, 
it is invariant under the action of $\Uh$, the algebra $A_t=A_{t,0}$
is a quantization of KKS bracket, and $A_{t,h}$ is a quotient
of $(S\g)_{t,h}$ by some ideal, $(S\g)_{t,h}\to A_{t,h}\to 0$.

As follows from \cite{DS}, there exists such an algebra
$A_{t,h}$ in case $M$ is a minimal orbit, i.e. $M$ is a hermitian
symmetric space. That this $A_{t,h}$ can be presented as a quotient
of $(S\g)_{t,h}$ follows easily from the view of irreducible components
of $A$. 

We are going to show here, that the problem also has a positive solution
for $M$ being a maximal orbit, i.e. can be defined as a set of zeros
of invariant functions in $S\g$. Such orbits are the orbits of diagonal
matrices with distinct elements on the diagonal,

The construction of $A_{t,h}$ is the following. There exists an isomorphism
of $\Uh$-modules $(S\g)_{h}\to W_h$, where $W_h=\oplus_kW_h^k$, the direct sum 
of the $k$-th symmetric parts of $T(\g_h)$ (see previous Section).
Consider the composition $W_h[t]\to T(\g_h)[t]\to(S\g)_{t,h}$. It is
an isomorphism, since it is an isomorphism at the point $h=0$.
It follows that $(S\g)_{t,h}$ is isomorphic to $W_h[t]$ as a $\Uh$-module,

Denote by $\I_{t,h}$ the submodule
of $\Uh$-invariant elements in $(S\g)_{t,h}$. It is obvious that $\I_{t,h}$
is isomorphic to $\oplus_k\I^k_h[t]$, where $\I^k_h$ is the invariant submodule
in $W^k_h$. Hence, $\I_{t,h}$ is a direct free $\C[[h]][t]$-submodule
in $(S\g)_{t,h}$. Moreover, $\I_{t,h}$ is a central subalgebra in 
$(S\g)_{t,h}$. Indeed, for $t\neq 0$ the algebra $(S\g)_{t,h}$ can be 
invariantly embedded in $\Uh$, but $\ad(\Uh)$-invariant elements in $\Uh$
form the center of $\Uh$. Yet $\I_{t,h}$ as an algebra
is isomorphic to $\I[[h]][t]$ with the trivial action of $\Uh$, where
$\I=\I_{0,0}$, the algebra of invariant elements in $S\g$ which is a 
polynomial algebra \cite{Dix} and, therefore, admits only the 
trivial commutative deformation. 

By the Kostant theorem \cite{Dix} $U(\g)$ is a free module over its center.
It follows that at the point $h=0$ the module $(S\g)_{t,0}$ 
is a free module
over the algebra $\I_{t,0}$. One can easily derive from this that
$(S\g)_{t,h}$ is a free module over $\I_{t,h}$.

Now, let the orbit $M$ be defined by invariant elements from $\I$.
Consider a character defined by $M$, the algebra homomorphism 
$\lambda:\I\to \C$ 
which sends each element from $\I$ to its value on $M$. 
Then, $\C$ may be considered as an $\I$-module, and the function algebra $A$ 
on $M$ is equal to $S\g/\Ker(\lambda)S\g=S\g\ot_\I\C$.
Extend the character $\lambda$ up to a character 
$\lambda_{t,h}:\I_{h,t}\to\C[[h]][t]$ in the trivial way and consider
$\C[[h]][t]$ as a $\I_{h,t}$-module. The tensor product over $\I_{t,h}$  
\be{*}
A_{t,h}=(S\g)_{t,h}\ot\C[[h]][t]
\ee{*}
is a $\C[[h]][t]$-algebra. It is a free $\C[[h]][t]$-module, since
$(S\g)_{t,h}$ is a free one over $\I_{t,h}$.  

It is obvious, that $A_{0,0}=A$, $A_{t,0}$ gives a quantization of
the KKS bracket on $M$,
and  $A_{t,h}$ is a quotient algebra of $(S\g)_{t,h}$.

\end{document}